\documentclass{amsart}

\usepackage{amsmath,amsthm,amssymb,amsfonts}

\newcommand{\R}{\mathbb{R}}

\newcommand{\un}{\mathbf{1}\!\!{\rm I}}
 
\newcommand{\be}{\begin{equation}}
\newcommand{\ee}{\end{equation}}
\newcommand{\bea}{\begin{eqnarray}}
\newcommand{\eea}{\end{eqnarray}}
\newcommand{\bean}{\begin{eqnarray*}}
\newcommand{\eean}{\end{eqnarray*}}
\newcommand{\rf}[1]{(\ref {#1})}

\def\xn{|\!|\!|} 
\def\La{\Lambda} 
\def\vt{\vartheta}
\def\dx{\,{\rm d}x}
\def\dy{\,{\rm d}y}
\def\dt{\,{\rm d}t}

\def\ds{\,{\rm d}s}

\def\e{\varepsilon}

\def\r{\varrho}
\def\th{\vartheta} 
\def\xn{|\!|\!|}

\newtheorem{theorem}{Theorem}

\newtheorem{lemma}[theorem]{Lemma}
\newtheorem{corollary}[theorem]{Corollary}

\theoremstyle{definition}

\def\proof{\noindent{\bf Proof\ \ }}
\def\qed{\hfill $\square$}

\theoremstyle{remark}

\def\int{\intop\limits}

\numberwithin{equation}{section}
\numberwithin{theorem}{section}

\author{Piotr Biler}
\address{\small Instytut Matematyczny, Uniwersytet Wroc\l awski,
 pl. Grunwaldzki 2/4, 50-384 Wroc\-\l aw, Poland}
\email{Piotr.Biler@math.uni.wroc.pl}

\author{Jacek Zienkiewicz}
\address{\small
 Instytut Matematyczny, Uniwersytet Wroc\l awski,
 pl. Grunwaldzki 2/4, 50-384 Wroc\-\l aw, Poland}
\email{Jacek.Zienkiewicz@math.uni.wroc.pl}

\title[existence for chemotaxis]{Existence of solutions for the Keller-Segel model of chemotaxis with measures as initial data}

\begin{document}

\begin{abstract}
A simple proof of  the existence of solutions for the two-dimensional Keller-Segel model with measures with all the atoms less than $8\pi$ as the initial data is given. This result has been obtained in \cite{SS-JFA} and \cite{BM} using different arguments. Moreover, we show a uniform bound for the existence time of solutions as well as an optimal hypercontractivity estimate. 
\end{abstract}



\date{\today}

\thanks{ The preparation of this paper was supported by the NCN grants   2013/09/B/ST1/04412, and DEC-2012/05/B/ST1/00692. 
\newline
Appeared in: Bull. Pol. Acad. Sci. Math. {\bf 63} (2015), 41--51. DOI: 10.4064/ba63-1-6.}

\maketitle

\baselineskip=16pt

\section{Introduction}
We consider in this paper the classical parabolic-elliptic Keller--Segel model of chemotaxis in two space dimensions
\bea
u_t-\Delta u+\nabla\cdot(u\nabla v)&=&0,\label{equ}\\
\Delta v+u &=& 0,\label{eqv}
\eea
supplemented with a nonnegative initial condition
\be
u(x,0)=u_0(x)\ge 0\label{ini}.
\ee
Here for $(x,t)\in {\mathbb R}^2\times[0,T)$, the function $u=u(x,t)\ge 0$ denotes the density of the population of microorganisms, $v=v(x,t)$ -- the density of the chemical secreted by themselves that attracts them and makes them to aggregate. The system \rf{equ}--\rf{eqv} is also used in modelling the gravitational attraction of particles in the mean field astrophysical models, see \cite{B-SM}.%, \cite{B-CM}. 

As it is well known,  cf. e.g. \cite{BDP},  the total mass of the initial condition
\be
M=\int u_0(x)\dx\label{M}
\ee
is the critical quantity for the global-in-time existence of nonnegative solutions.  
Namely, if $M\le 8\pi$, then solutions of \rf{equ}--\rf{ini} (with $u_0$ --- a~finite nonnegative measure)
 exist for all $t\ge 0$. For the local-in-time existence, it should be assumed that all the  atoms of the measure $u_0$ are of mass less than $8\pi$, see \cite[Th. 2]{BM}. When $M>8\pi$, nonnegative solutions blow up in a finite time, and  for radially symmetric solutions  mass equal to $8\pi$ concentrates at the origin at the blowup time, see e.g. \cite{BKZ}. 

Our goal in this note is to give an alternative proof of the local in time existence of solutions to \rf{equ}--\rf{ini} when $u_0\in{\mathcal M}(\R^2)$ is a nonnegative finite measure with all its atoms of mass less than $8\pi$. 
We believe that this approach is conceptually simpler than that in the recent paper \cite{BM} (which used elaborated arguments for interactions of solutions emanating from localized pieces of initial data), and those in previous papers \cite{SS-ADE}, \cite{SS-JFA}. 
The latter approaches used heavily the free energy functional for system \rf{equ}--\rf{eqv} considered in bounded domains. 
Moreover, our condition \rf{small} seems to be more clear, and shows that measures with small atoms, which are not well separated as it was assumed in \cite{BM}, are also admissible as initial data for the system \rf{equ}--\rf{eqv}. Compared to \cite{BM} here, however, we   obtain neither the uniqueness property of solutions nor the Lipschitz property of the solution map. 

The main result of this paper is 

\begin{theorem}\label{main}
Let $0\le u_0\in L^1(\R^2)\cap L^\infty(\R^2)$ be a smooth initial density for \rf{equ}--\rf{eqv} such that 
\be
\|u_0\ast \un_{B(1)}\|_\infty\le 8\pi-\e_0\label{small}
\ee 
for some fixed 
$\e_0>0$ and the unit ball $B(1)$ centered at the origin of  
in $\R^2$. 
Then, there exists a solution of the problem \rf{equ}--\rf{ini} on the interval $[0,t_0]$ with $t_0=t_0(\e_0,M)$  such that 
\be
\sup_{0<t\le t_0}t^{1-1/p}\|u(t)\|_p\le B,\label{hiper}
\ee
where the constant $B$ depends on $M$ and $\e_0$ (in particular, $B$ does not depend on $\|u_0\|_\infty$). 
\end{theorem}

Note that the condition \rf{small} reads 
$$
\int_{B(x,1)}u_0(y)\dy=\int_{B(1)}u_0(x-y)\dy \le 8\pi-\e_0
$$
for all the   balls $B(x,1)$ of radius $1$ centered at an arbitrary $x\in\R^2$ and this, in particular, means that if, more generally, $u_0$ were a nonnegative measure then its atoms will be strictly less than $8\pi$. 
The key property of our estimate of $t_0$ is that it depends only on $M$ and $\e_0$   for all $u_0$ satisfying \rf{small} with a given $M$ in \rf{M}.%\rf{mass}. 

We recall, that for each $\lambda>0$ and each solution $u$ of \rf{equ}--\rf{eqv} of mass $M$ the function 
\be
u_\lambda(x,t)=\lambda^2u(\lambda x,\lambda^2t)\label{scale}
\ee 
is also a solution, with its mass again equal to $M$. 

Of course, by a suitable  scaling \rf{scale} of initial data we see that we may satisfy the assumptions of the result on the local existence in Theorem \ref{main} for  any nonnegative $u_0\in{\mathcal M}(\R^2)$ with its atoms strictly less than $8\pi$. Then, it is clear that we arrive at the following corollary, cf. \cite[Theorem 2]{BM}. 

\begin{corollary}
The system \rf{equ}--\rf{eqv} has a local-in-time solution for each initial nonnegative finite measure $u_0$ with all its atoms strictly less than $8\pi$.
\end{corollary}

Indeed, it is sufficient to approximate (in the sense of the weak convergence of measures) such a measure $u_0$ by a sequence of initial data satisfying (after the rescaling \rf{scale} with a single $\lambda>0$) all the assumptions of Theorem \ref{main}. 
This approximation is possible by taking, e.g., ${\rm e}^{\delta_n\Delta}u_0$ for any sequence $\delta_n\searrow 0$. 
Then, the existence time $t_0$ is bounded from below by a positive quantity (since $t_0$ depends on $M$ and $\lambda$ only). Next, we infer from the hypercontractivity estimate \rf{hiper} and from the standard regularity theory for parabolic equations that for every multiindex $\alpha$ 
$$
\|D^\alpha u(t)\|_p\le C_\alpha B\,t^{1/p-1-|\alpha|/2},
$$ 
which permits us to pass to the limit with (a subsequence of) the approximating solutions which are, in fact, smooth on $\R^2\times(0,t_0)$. We obtain in such a way a  solution to \rf{equ}--\rf{eqv} with the measure $u_0$, and this solution is also smooth on $\R^2\times(0,t_0)$.

\bigskip

The proof of Theorem \ref{main} will be a consequence of a well-known fact  in \cite{B-SM,BCKZ} on the estimate of the existence time for a solution by mass of the initial condition only, see \ref{stare} in Sec. 2,  by using  a rather delicate argument of localization repeatedly. 

The existence of solutions results are proved (e.g. as in \cite{B-SM}) for the integral formulation of the system \rf{equ}--\rf{ini}
\be
u(t)={\rm e}^{t\Delta}u_0+B(u,u)(t),\label{D}
\ee
whose solutions are called {\em mild} solutions of the original Cauchy problem. 
Here, the bilinear term $B$ is defined as 
\be
B(u,z)(t)=-\int_0^t\left(\nabla{\rm e}^{(t-s)\Delta}\right)\cdot\left(u(s)\, \nabla(-\Delta )^{-1}z(s)\right)\ds.\label{form}
\ee 
It is well known that the heat semigroup ${\rm e}^{t\Delta}$, satisfies the following $L^q-L^p$ estimates 
\be
\|{\rm e}^{t\Delta}z\|_p\le Ct^{1/p-1/q}\|z\|_q\label{lin1}
\ee
and 
\be 
\|\nabla{\rm e}^{t\Delta}z\|_p\le Ct^{-1/2+1/p-1/q}\|z\|_q\label{lin2}
\ee
 for all $1\le q\le p\le\infty$. 
 Moreover, for each $p>1$ and $z\in L^1(\R^2)$ the following relation holds 
 \be
 \lim_{t\to 0}t^{1-1/p}\|{\rm e}^{t\Delta}z\|_p=0.\label{function}
 \ee
 This is the consequence of, e.g., much more general inequality valid for  every finite measure $\mu\in{\mathcal M}({\mathbb R}^2)$ and every $p>1$  
\be
\limsup_{t\rightarrow 0}t^{1-1/p}\|{\rm e}^{t\Delta}\mu\|_p\le C_p\|\mu_{\rm at}\|_{{\mathcal M}(\R^2)}\equiv C_p\sum_{\{x: \mu(\{x\})\ne 0\}}|\mu(\{x\})|, \label{measure}
\ee
where $\mu_{\rm at}$ denotes the purely atomic part of the measure $\mu$. The proof of \eqref{measure} is contained in \cite[Lemma 4.4]{GMO}.
This fact, equivalent to the condition \rf{small} rescaled to other balls of a fixed radius, see \rf{small2} below,   is crucial in the analysis of applicability of the Banach contraction argument to the equation \rf{D}. 
 
We recall that the formulation of our existence results %Theorem \ref{stare} 
in \cite{B-SM} used in fact condition \rf{measure} in the definition of the functional space where solutions have been looked for: $\{ u:(0,T)\to L^p(\R^2):\ \xn u\xn\equiv \sup_{0<t<T}t^{1-1/p}\|u(t)\|_p<\infty\}$, and then a smallness condition has been assumed on the quantity $\xn {\rm e}^{t\Delta}u_0\xn$. 
\medskip  

The heuristics behind the argument leading to the proof of Theorem \ref{main} is the following: 
the initial data diffuse into a domain which size grows as $t^{1/2}$ in time as in Corollary \ref{disp}. 
Thus, we need to find a moment of time $\tau\ge 0$  when a counterpart of the condition \rf{small}
\be
\|u(\tau)\ast\un_{B(\r)}\|\le m_0  
\label{small2}
\ee 
 holds with a sufficiently small $m_0$ given in \rf{ast1} and $\r>0$ suitably small in order to apply the local existence result in Theorem \ref{stare}.  
\bigskip

\noindent 
{\bf Remark.} 
When  equation \rf{eqv} is replaced by the nonhomogeneous heat equation $\tau v_t=\Delta v+u$ (and thus we consider the parabolic-parabolic version of the Keller-Segel model), the situation seems be more complicated. 
For instance, if $\tau\gg 1$, then there exist global-in-time solutions with $M>8\pi$ which emanate from $M\delta_0$ as (purely atomic) initial data. These are self-similar solutions which are regular and {\em nonunique} for sufficiently large $M$, cf. e.g. \cite{BCD} and comments in \cite{BGK}. 
\bigskip

We will apply in the proof of Theorem \ref{main}  simple (but rather subtle) techniques of weight functions and scalings. The core of our analysis are the uniform (with respect to the initial distributions) estimates on the maximal existence time, expressed in  terms of dispersion of the initial data. 

\bigskip

\noindent
{\bf Notations.}
The integrals with no integration limits are understood as over the whole space $\R^2$: \ \ $\int=\int_{\R^2}$. The letter $C$ denotes various constants which may vary from line to line but they are independent of solutions. 
The norm in $L^p(\R^2)$ is denoted by $\|\, .\, \|_p$. 
 The kernel of the  heat semigroup on $\R^2$,  denoted by ${\rm e}^{t\Delta}$,  is given by  $G(x,t)=(4\pi t)^{-1}\exp\left(-\frac{|x|^2}{4t}\right)$. 
\medskip

\section{Proof of Theorem \ref{main}} 

The proof of the estimate of the existence time in Theorem \ref{main} is split into several lemmata.

 For any fixed $x_0\in\R^2$ we define {\em the  local moment} of a solution $u$ by 
\be
\La(t) \equiv  \int\psi(x-x_0)u(x,t)\dx. \label{mom}
\ee
Here the weight function
\be
\psi(x)=\left(1-|x|^2\right)^2_+\ \ {\rm with}\ \ \nabla\psi(x)=-4x(1-|x|^2)_+\ \ {\rm and}\ \ \Delta\psi(x)=16|x|^2-8\ge -8, \label{defpsi} 
\ee 
is a fixed radial, piecewise ${\mathcal C}^2$, nonnegative function $\psi$, supported on the unit ball such that $\psi(0)=1$.
Our particular choice of the function $\psi$ is not critical.

\begin{lemma}\label{L0}
Suppose that $w=w(x)$ is a nonnegative function  locally in $L^1\cap L^\infty$, 
$
\int_{B(1)} w(x)\dx\le m,
$ 
and $\r,\, \delta\in(0,1)$. 

\noindent 
(i) Then there exists a number $H_0 \in(0,1)$ such that  $\int_{B(\r)}w(x)\dx\le (1-\delta)m$ implies $\int \psi(x)w(x)\dx\le (1-H_0)m$. 

\noindent
(ii) Similarly, there exists $H_1\in(0,1)$ such that if $\int_{B(1)} w(x)\dx\le m$ and  $\int \psi(x)w(x)\dx\ge (1-H_1)m$, then $\int_{B(\r)}w(x)\dx \ge 
 \left(1-\frac\delta{2}\right)m$. 

\noindent 
(iii) Suppose that the inequality $\int \psi(x)w(x)\dx\le (1-H)m$ holds with some $H\in(0,1)$. Then the bound $\int_{B(\beta)}w(x)\dx\le\left(1-\frac{H}{2}\right)m$ holds for $\beta^2\le \frac{H}{4}\le \frac14$. 
\end{lemma}

\proof 
The properties (i)--(iii) are simple consequences of \rf{defpsi}. 
Indeed, 
\bea
\int \psi(x)w(x)\dx &\le& \int_{B(\r)}w(x)\dx+\sup_{B(1)\setminus B(\r)}\psi(x)\cdot \int_{B(1)\setminus B(\r)} w(x)\dx\nonumber\\
&\le& \int_{B(\r)}w(x)\dx+(1-\r^2)^2\int_{B(1)\setminus B(\r)} w(x)\dx \nonumber\\
&=&(1-\r^2)^2\int_{B(1)}w(x)\dx+\left(1-(1-\r^2)^2\right)\int_{B(\r)}w(x)\dx\nonumber\\ 
&\le& (1-\r^2)^2m +\left(1-(1-\r^2)^2\right)(1-\delta)m\nonumber\\
&=&(1-H_0)m,\nonumber
\eea
where $1-H_0=(1-\r^2)^2+\left(1-(1-\r^2)^2\right)(1-\delta)=1-\delta \left(1-(1-\r^2)^2\right)$. 
\medskip

(ii) is equivalent to (i) with $\delta$ replaced by $\frac12\delta$. 
\medskip

(iii) For $|x|\le \beta$, $\beta^2\le\frac{H}{4}$ and $H\le 1$ the inequality $\psi(x)\ge \frac{1-H}{1-\frac{H}{2}}$ holds. 
\qed
\medskip

Next,  we show a result on the dispersion of the initial data evolving according to \rf{equ}--\rf{eqv}. 

\begin{lemma}\label{L1}
Let $u$ be a solution to \rf{equ}--\rf{eqv} such that for $t\in[0,A]$ 
\be
\|u_0\ast \un_{B(R_0)}\|_\infty\le m\label{small-m}
\ee 
for some $A>0$, $R_0=6\cdot128\frac{\pi M}{\e_0}>0$ and $m_0\le m\le 8\pi-\e_0$. 
Then, there exist numbers $A_1=A_1(M,\e_0)$, $\delta=\delta(M,\e_0,m_0)$ and $\varrho=\r(M,\e_0,m_0)$ such that if $\int_{|y-x_0|\le\r}u(y,t)\dy\ge (1-\delta)m$ for some $t\in [0,A]$, then the differential  inequality 
$$\La'(t)\le -\vartheta$$ 
holds with some $\vt=\vt(M,\e_0,m_0)>0$. 
\end{lemma}

\proof 
First we  give a uniform estimate of  the time derivative of the moment $\La(t)$: 
$$
|\La'(t)|\le C_M. 
$$ 
Let us compute the time derivative of $\La$ using equations \rf{equ}--\rf{eqv} and \rf{defpsi}. 
Symmetrizing the bilinear integral $\int u(x,t)\nabla v(x,t)\cdot\nabla\psi(x)\dx$ with the solution $v$ of \rf{eqv} given by \ \  $v(x,t)=-\frac{1}{2\pi}\int u(y,t)\log |x-y|\dy$ we obtain 
\bea
\La'(t)&=&\int u(x,t)\Delta\psi(x)\dx \nonumber\\
&\quad&+\frac{1}{4\pi}\iint \frac{\nabla\psi(x)-\nabla\psi(y)}{|x-y|^2}\cdot(x-y)u(x,t)u(y,t)\dx\dy. \label{lambda-ev}
\eea
From \rf{lambda-ev} and \rf{defpsi} we immediately get $\La'(t)\le 8M+4M^2$. 
\medskip

 Using \rf{defpsi},  the bound $|\nabla\psi(x)-\nabla\psi(y)|\le 4$ and the relation $\int_{B(\r)}u(x,t)\dx\le 8\pi$ with $\r\le 1<2\le R_0$ and $\frac{1}{R_0-\r}\le \frac{2}{R_0}$, we arrive at 
\be
\left|\int_{|x|<\r}\int_{2\le R_0<|y|}\frac{4x(1-|x|^2)_+
-4y(1-|y|^2)_+}{|x-y|^2}\cdot(x-y)u(x,t)u(y,t)\dx\dy\right|\le \frac{16\pi M}{R_0}\cdot 4.\label{i} 
\ee
\medskip

 Next, in the annulus $\r<|y|\le R_0$ we have 
\be
\left|\int_{|x|<\r}\int_{\r<|y|\le R_0}\frac{\nabla\psi(x)-\nabla\psi(y)}{|x-y|^2} \cdot(x-y) u(x,t)u(y,t)\dx\dy\right|\le B\, 8\pi\, \delta, \label{ii}
\ee
where we applied the bound 
\be 
|\nabla\psi(x)-\nabla\psi(y)|\le B|x-y|\label{nablapsi}
\ee  
which holds for some constant $B$, as well as $\int_{\r<|y|\le R_0}u(y,t)\dy\le \delta m<\delta\, 8\pi$. 
\medskip

  Finally, by \rf{lambda-ev} we have simply 
$$
\left|\int u(x,t)\Delta\psi(x)\dx+8\int_{B(\r)}u(x,t)\dx\right|\le 64\pi\delta+16\r^2\, 8\pi.
$$

Now, the crucial estimate for the bilinear integral in \rf{lambda-ev} is 
\bea
\left|\int_{|x|<\r}\int_{|y|<\r} \frac{4x(1-|x|^2)_+-4y(1-|y|^2)_+}{|x-y|^2}\cdot(x-y)u(x,t)u(y,t)\dx\dy \right.\nonumber\\
\quad \left.-4 \left(\int_{B(\r)}u(x,t\dx \right)^2 \right|\le B\delta M^2.
\eea
Here we used the following properties of the weight function $\psi$:
$$
|\psi(x)-1|=\left|(1-|x|^2)_+^2-1\right|=|2|x|^2-|x|^4|_+\le 2|x|^2,\ \ |\Delta\psi(x)+8|\le 16|x|^2,$$
and an improvement of \rf{nablapsi}: 
$$
|\nabla\psi(x)-\nabla\psi(y)+4(x-y)|\le B\r|x-y|
$$
valid for all $|x|,\, |y|\le\r$. 
Therefore, we get 
$$
\La'(t)\le -8\int_{B(\r)}u(x,t)\dx+\frac1\pi\left(\int_{B(\r)}u(x,t)\dx\right)^2 +\frac{128\pi M}{R_0} +16\pi B\delta +64\pi\delta+16\r^28\pi.
$$
Since $$\frac1\pi\int_{B(\r)}u(x,t)\dx\left(-8\pi+\int_{B(\r)}u(x,t)\dx\right)\le -\e_0,$$ 
it suffices to choose $R_0=\frac{6C_0}{\e_0}$, $\delta\le \frac{\e_0}{6C_1}$, $\r^2\le \frac{\e_0}{6C_2}$, therefore  $R_0=\frac{6\cdot 128\pi M}{\e_0}$. 
\qed
\bigskip

Note that a variant of Lemma \ref{L3} below has been obtained in \cite{SS-JFA} by a  (rather elaborate) radial rearrangement argument of \cite{DNR}. 
The proof we present uses only weight functions and localized moments defined by them.

\begin{lemma}\label{L3}
Suppose that $u=u(x,t)$ is a solution of \rf{equ}--\rf{eqv} satisfying all the assumptions of Lemma \ref{L1}. Then for all $t\ge A_0= A_0(M,m_0,\e_0)$ and $H_1=H_1(M, m_0,\e_0)$ we have 
$$
\int_{B(\r)}u(x,t)\dx \le (1-H_1)m.
$$ 
\end{lemma}

\proof 
Let $\tau_0=0$. If $\int_{B(\r)} u(x,\tau_0)\dx\le (1-\delta)m$, then $\tau_1=\tau_0$. Otherwise, let 
$$
\tau_1=\inf\left\{ \tau<A: \int_{B(\r)}u(x,\tau)\dx=(1-\delta)m\right\}.
$$
In order to obtain necessary estimates for $\tau_1$, observe that 
by Lemma \ref{L1}, for all $t\in[\tau_0,\tau_1]$ the inequality 
$\frac{\rm d}{\dt}\La(t)\le-\vt$ holds. Since $\La(0)\le m$ and $\La(t)\ge 0$, we have $\tau_1-\tau_0\le m/\vt$. By the first part of Lemma \ref{L0} we arrive at the inequality $\La(\tau_1)\le (1-H_0)m$. 
Next, we define either 
$$\tau_2=\inf \left\{\tau_1<\tau<A: \La(\tau)=(1-H_1)m\right\}
$$
if this exists, or $\tau_2=A$.  
Then for every $t\in[\tau_1,\tau_2]$ we obtain $\La\le(1-H_1)m$. 
If $\tau_2=A$, we are done. If not, 
by the second part of Lemma \ref{L0} we infer that 
$\int_{B(\r)} u(x,\tau_2)\dx\ge \left(1-\frac\delta{2}\right)m>(1-\delta)m$. 
Therefore, $\La(\tau_2)\le-\vt$ implies that $\La(\tau_2-h)>(1-H_1)m$ for a sufficiently small $h>0$, contrarily to the definition of $\tau_2$. Thus we get $\La(t)\le (1-H_1)m$ for $t\in [\tau_1,A]$. 
Consequently, by Lemma \ref{L0}, we have $\int_{B(\beta\r)}u(x,t)\dx\le \left(1-\frac{H_1}{2}\right)m,$ for $t\ge \frac{m}{\th}\equiv A_0(M,\e_0)$ and $\beta<\frac12H^{1/2}$ as in Lemma \ref{L0} (iii). We denote $\beta\r$ in the sequel again by $\r$. 
\qed
\bigskip

\begin{corollary}\label{cor1}
If $u$ solves the system \rf{equ}--\rf{eqv} on the time interval $[0,AT]$, $A\ge A_0$,  and satisfies the estimate 
\be
\int_{B(T^{1/2})}u(x,t)\dx\le m<8\pi-\e_0,\label{par}
\ee
then for $t\in[A_0T,AT]$ we have $\int_{B(\r T^{1/2})}u(x,t)\dx\le (1-H_1)m$. 
\end{corollary}

\proof 
This follows from Lemma \ref{L3} applied to the rescaled function $u_T(x,t)=Tu(T^{1/2}x,Tt)$ which, evidently,  is also a solution of \rf{equ}--\rf{eqv}.  \qed
\bigskip

\begin{corollary}\label{cor2}
If a solution $u$ of \rf{equ}--\rf{eqv} exists on the time interval $[0,A]$ then for all $t\in[A_0(1+\r_1+\dots+\r_1^{j-1}),A)$  with $\r_1=\frac{\r}{R_0}$ the function $u$ satisfies the estimate 
$$
\int_{B(\r^j_1)}u(x,t)\dx\le (1-H_1)^jm
$$ 
as long as the inequality $(1-H_1)^jm\ge m_0$ holds. 
\end{corollary}

\proof 
It suffices to apply Corollary \ref{cor1}  to the functions $u(x,t)$, $u(x,A_0+t)$, $u(x,A_0+\r_1 A_0+t)$, $\dots$, rescaled with $T_0=1$, $T_1=\r_1^2$, $T_3=\r_1^4$, $T_3=\r_1^6$, $\dots$,  consecutively. \qed
\bigskip

\begin{corollary}\label{cor3}
Suppose that a solution $u$ of the system \rf{equ}--\rf{eqv} exists for $t\in[0,A]$ and satisfies 
$$
\|u(t)\ast\un_{B(1)}\|_\infty\le 8\pi-\e_0.
$$ 
Then for any $s\in[0,1]$ and $t\in\left[\frac{s^2}{1-\r_1}A_0,A\right]$ the inequality 
$$
\|u(t)\ast\un_{B(s r)}\|_\infty \le m_0,
$$ 
holds with $r=\r_1^{\left[\left|\frac{\log 8\pi}{\log(1-H_1)}\right|\right]+1}$. 
\end{corollary}

\proof 
Apply Corollary \ref{cor2} to the rescaled solution.  \qed

\bigskip

Now we recall the existence  result in \cite{B-SM,BCKZ}

\begin{theorem}\label{stare} 
There exists a  (small) $m_0>0$ such that the condition 
\be
\|u_0\ast\un_{B(1)}\|_\infty\le 2m_0\label{ast1}
\ee 
guarantees the existence of a local-in-time solution (on a time interval $[0,T]$ with $T=T(M)$) satisfying the %decay 
estimate $\|u(t)\|_p\le Ct^{1/p-1}$ for $p=\frac43$.
\end{theorem} 

In fact, we will use in the sequel the following immediate  version of Theorem \ref{stare} which takes into account the scale invariance \rf{scale}  

\begin{corollary}\label{disp}
For each $\sigma>0$ there exists $\alpha>0$ such that the condition 
$$\|u_0\ast\un_{B(\sigma\tau^{1/2})}\|_\infty\le 2m_0$$ 
implies the existence of the solution of  \rf{equ}--\rf{ini} on the time interval $[0,\alpha\tau]$. Here $\tau>0$ is any small positive number and $\alpha$ can be chosen as  $\alpha=\alpha_0\sigma^2$ for some $\alpha_0>0$.  
\end{corollary}

Theorem \ref{stare} has  been proved (even for sign-changing measures) in \cite[Theorem 2]{B-SM} (cf. also \cite[proof of Theorem 2.2]{BCKZ}) using a standard contraction argument applied to the formulation \rf{D}. 

\begin{lemma}\label{exi}
Suppose that  $u_0\in L^1(\R^2)\cap L^\infty(\R^2)$ is a smooth nonnegative function satisfying the condition 
$$
\|u_0\ast \psi\|_\infty \le 8\pi-\e_0.
$$
Then, the solution $u$ with the initial condition \rf{ini}  $u_0$ exists at least on the time interval $\left[0,\frac{\e_0}{2C_M}\right]$.        
\end{lemma}

\bigskip 

\proof \ \ 
The inequality $\La'(t)\le C_M$ implies that $\La(t)\le 8\pi-\e_0/2$ for $t\le \tau_1\equiv \min\left\{\frac{\e_0}{2C_M}, \tau\right\}$, where $\tau$ is the maximal existence time of $u$. 
By  (iii) of Lemma \ref{L0} we obtain the bound 
$\|u(t)\ast\un_{B(\r_1)}\|_\infty\le 8\pi-\frac{\e_0}{4}$ for all $t\le \tau_1$. 
From Corollary \ref{cor3} we infer that there exists  $\sigma_0=\eta\tau_1^{1/2}$ such that $\|u(\tau_2)\ast \un_{B(\sigma_0)}\|_\infty\le m_0$ for $\frac{\tau_1}2<\tau_2<\tau_1$.  
By Theorem \ref{stare} the solution with $u_0=u(x,\tau_2)$ as the initial condition exists for $t\in[0,\alpha\tau_2]$ with some $\alpha>0$ independent of $\tau_2$. 
Therefore, by Theorem \ref{stare}, this solution can be continued onto the interval $[0,\tau_1+\alpha\tau_2]$. This solution satisfies the estimate $\|u(\tau_1+\alpha\tau_2)\|_p\le C\tau_2^{1/p-1}$ for each $p\in\left[\frac43,2\right]$. 
Finally, a recurrence argument permits us to obtain a classical solution $u=u(x,t)$  on the whole interval  $[0,T_0]$ with $T_0=T_0(\e_0,M)$, and applying once more Corollaries \ref{cor2} and \ref{disp} this satisfies the hypercontractive estimate 
$\|u(t)\|_p\le Ct^{1/p-1}$ for $p\in\left[\frac43,2\right]$. 
The extrapolation of that estimate to the whole range of $p\in(1,\infty)$ is standard, see e.g. \cite{BCKZ}. 
\qed

\end{document}